\documentclass[11pt,reqno]{amsart}

\usepackage[utf8]{inputenc}
\usepackage[T1]{fontenc}
\usepackage{newunicodechar}
\newunicodechar{–}{\textendash}
\usepackage{amsmath,amssymb,amsthm}
\usepackage{enumitem}
\usepackage[colorlinks=true,linkcolor=blue,citecolor=blue,urlcolor=blue]{hyperref}

\theoremstyle{plain}
\newtheorem{theorem}{Theorem}
\newtheorem{proposition}[theorem]{Proposition}
\newtheorem{lemma}[theorem]{Lemma}
\newtheorem{corollary}[theorem]{Corollary}
\theoremstyle{remark}
\newtheorem{remark}[theorem]{Remark}

\newcommand{\Z}{\mathbb{Z}}
\newcommand{\rot}{\operatorname{rot}}
\newcommand{\rev}{\operatorname{rev}}
\newcommand{\tr}{\operatorname{tr}}
\newcommand{\T}{^{\mathsf{T}}}
\begin{document}
\sloppy

\title{Determinants of amphichiral knots of braid index at most three}

\author{Werner Alois Stanggassinger}
\address{Wasserburg am Inn, Bavaria, Germany}
\email{stanggassinger@stainlesssteel4u.de}

\thanks{This work was developed in collaboration with the AI system Claude Code (Anthropic):
the choice of problem, the direction of the work and the computing hardware are due to the author
(a sculptor); the formalisation, the algorithms and the proofs were developed by Claude Code
(Anthropic) working with the author. The collaboration is stated as a fact, not a disclaimer; details
are in the section on tool use at the end.}


\date{12 November 2026. Not refereed.}

\begin{abstract}
We show that an odd number $n$ is the determinant of an amphichiral knot of braid index at most three if and only if both $n$ and $n+4$ are sums of
two squares. The proof rests on the known description of amphichiral knots of braid index three as closures of alternating 3-braid words with a
rotation or reflection symmetry (Birman--Menasco, Stoimenov, Jung). Determinant formulas for closed 3-braids are also known (e.g.\ through the reduced
Burau representation at $t=-1$); what is added here is their evaluation on the two symmetric families, which gives the condition on $n+4$, and an
explicit construction for the converse. For amphichiral knots of braid index three this describes which determinants occur, which bears on the determinant part of a question of Jung. We did not find the statement in the literature we searched. As a consequence, for an odd $n$ with $n+4$
not a sum of two squares, every amphichiral knot with determinant $n$ has braid index at least five; we test this on 4,999 amphichiral knots of
Burton's census with up to 19 crossings whose determinant satisfies this hypothesis.
\end{abstract}

\maketitle

\section{Statement}

A knot $K$ is $+$amphichiral if it is isotopic to its mirror image, $-$amphichiral if it is isotopic to the reverse of its mirror image, and
amphichiral if one of the two holds. The determinant of a knot is odd, and for every amphichiral knot it is a sum of two squares~\cite[Remark~(7), via
Goeritz]{HK}.

\begin{theorem}\label{thm:main}
$\{\det K : K \text{ an amphichiral knot of braid index} \le 3\} = \{\,n \text{ odd} : n \text{ and } n+4 \text{ are sums of two squares}\,\}$.
\end{theorem}

Here $0$ counts as a square (so $5$, the determinant of the figure-eight knot, belongs to the right-hand side: $5=2^2+1^2$, $9=0^2+3^2$).
Among the odd sums of two squares the first ones excluded by the condition on $n+4$ are $17$, $29$, $53$, $65$, $73$.
Braid index~$1$ and~$2$ contribute only the unknot ($n=1$), since the torus knots $T(2,p)$, $|p|\ge3$, are chiral.

\begin{corollary}\label{cor:b5}
If $n$ is odd and $n+4$ is not a sum of two squares, every amphichiral knot with determinant $n$ has braid index at least $5$.
\end{corollary}

\begin{proof}
By Theorem~\ref{thm:main} the braid index is at least $4$, and the braid index of an amphichiral knot is odd~\cite[Theorem~A]{J} (the proof is one line: by the theorem of Dynnikov and Prasolov~\cite[Theorem~9]{DP} (the Jones conjecture) all minimal braids of a
knot have the same exponent sum, so a minimal braid of an amphichiral knot has exponent sum $0$~\cite[Lemma~3.3 and eq.~(1)]{J}, and the parity of the
permutation of a braid whose closure is a knot then makes the number of strands odd).
\end{proof}

Jung asks~\cite[Question~2]{J} whether invariants such as the determinant can be ``detected or computed from the braid word''. For amphichiral knots of
braid index three, Lemmas~\ref{lem:N} and~\ref{lem:P} below compute the determinant from the word, and Theorem~\ref{thm:main} describes the values
that occur. (The author of~\cite{St15,StCAG,StMZ,St06,StW} is the same person; we keep the spelling of each publication.)

\section{Conventions and what is used from the literature}

$B_3=\langle\sigma_1,\sigma_2\rangle$, $\Delta=\sigma_1\sigma_2\sigma_1$, $C=\Delta^2=(\sigma_1\sigma_2)^3$. For a sequence $s=(e_0,\dots,e_{L-1})$ of
positive integers, $L$ even, the \emph{alternating word} of $s$ is $\sigma_1^{e_0}\sigma_2^{-e_1}\sigma_1^{e_2}\cdots\sigma_2^{-e_{L-1}}$. Let
$U(e)=\begin{pmatrix}1&e\\0&1\end{pmatrix}$, $L(e)=\begin{pmatrix}1&0\\e&1\end{pmatrix}$, $J=\begin{pmatrix}0&1\\1&0\end{pmatrix}$ and
$M(s)=U(e_0)L(e_1)U(e_2)\cdots L(e_{L-1})$. These are the matrices of the reduced Burau representation at $t=-1$ ($\sigma_1\mapsto U(1)$,
$\sigma_2\mapsto L(-1)$, so $\sigma_1^{e}\mapsto U(e)$ and $\sigma_2^{-e}\mapsto L(e)$). For a closed 3-braid knot with Burau matrix $M$ at $t=-1$ one has $\det=|2-\tr M|$: Datta~\cite[\S5]{D} states
$\Delta_L(-1)=i^{w}(\theta-2)$ for the closure of any 3-braid, $w$ its writhe and $\theta$ the trace of its matrix (with our convention), and it also follows
from $(1+t+t^2)\Delta(t)\doteq\det(I-\psi_t(\beta))$; we checked it against SnapPy on random 3-braid knots. For $M(s)=\begin{pmatrix}a&b\\c&d\end{pmatrix}$
all entries are nonnegative and $ad=1+bc\ge1$, so $\tr M(s)=a+d\ge2$ and $\det=\tr M(s)-2$. General determinant formulas for closed 3-braids are in~\cite[Theorem~3.1,
from Qazaqzeh--Chbili]{Ch}. $\rot_r(s)=(e_r,e_{r+1},\dots)$ (indices mod $L$) and $\rev(s)=(e_{L-1},\dots,e_0)$.

\emph{Types.} $s$ is of \emph{positive type} if $s=\rot_r(s)$ for some odd $r$, and of \emph{negative type} if $e_i=e_{k-i}$ for all $i$ for some odd~$k$
(equivalently, $s$ is $\rev(s)$ rotated by an even amount).

\emph{From the literature.} (i) Birman and Menasco classified closed 3-braids~\cite{BM}: every conjugacy class of $B_3$ contains a Schreier normal form,
unique up to cyclic permutation of the syllables after $C^k$, either $C^k\sigma_1^{-p_1}\sigma_2^{q_1}\cdots\sigma_1^{-p_r}\sigma_2^{q_r}$ ($r,p_i,q_i\ge1$,
``generic'') or one of $C^k\sigma_1^{p}$, $C^k\sigma_1\sigma_2$, $C^k\sigma_1\sigma_2\sigma_1$, $C^k\sigma_1\sigma_2\sigma_1\sigma_2$~\cite[pp.~104--105]{BM};
oriented links of braid index $\le3$ are determined by the conjugacy class of a 3-braid representative up to the unknot, the torus links of type
$(2,p)$ and the flype case, in which a word and its reverse represent the same link~\cite[pp.~76--77, 106--107]{BM}. Their Corollary~7.1~B~\cite[p.~107]{BM}
describes the $+$amphichiral generic closures by the condition that $\rot_1(A)$ agrees with the array $A$ up to an even rotation, i.e.\ $A=\rot_{\mathrm{odd}}(A)$ (positive type), and $k=0$; it has no flype clause --- $6_3$, the closure of
$\sigma_1^{-1}\sigma_2^2\sigma_1^{-2}\sigma_2$, is $+$amphichiral but only of negative type --- and does not treat $-$amphichirality.
(ii) \emph{Amphichiral 3-braid links.} Stoimenov~\cite[Theorem~1.1(a)]{St15} proves that an orientedly amphicheiral 3-braid link (positively or negatively~\cite[Definition~2.2]{St15}) is the closure of an
alternating 3-braid whose Schreier vector $S=(p_1,-q_1,\dots,p_s,-q_s)$ admits a ``dihedral (anti)symmetry'', i.e.\ can be turned into $-S$ by cyclic
permutation and inversion~\cite[p.~448]{St15}; a rotation by $r$ changes all signs iff $r$ is odd, and a reflection $i\mapsto k-i$ does iff $k$ is
odd, so these are exactly our two types. His proof bypasses the method of~\cite{BM} (his Remark~3.2; it goes back to~\cite[Theorem~6.2]{St06}), relying on~\cite{BM} only for where flypes can occur (for the equality of writhes he cites~\cite{BM} or, alternatively, his own earlier work). For prime
knots, Jung~\cite[Theorems~4.13 and B; the symmetry type in Theorem~C]{J} gives the same description, with the flype case $(k,1,1,k)$; in his terms positive type is
``odd period'' and negative type ``palindrome''~\cite[\S4, eq.~(2)]{J}. The two families, with the proof that their closures are amphichiral
(our Lemma~\ref{lem:S}), also appear (for words with at least two pairs of syllables) as family~(3) in Lisca's theorem on 3-braid knots of finite concordance order~\cite[Theorem~1.1]{Li}; he does not consider
which determinants occur. That the alternating half follows from~\cite{BM} is also noted
in~\cite[Prop.~7.1(2)]{StMZ}. For completeness we give the short derivation from~\cite{BM} in Proposition~\ref{prop:A}.

\section{Lemmas}

\begin{lemma}[symmetry]\label{lem:S}
Up to conjugation, the mirror image followed by conjugation with $\Delta$ sends the alternating word of $s$ to that of $\rot_1(s)$, and reading the
word backwards sends it to that of $\rot_1(\rev s)$. Rotation of $s$ by an even amount is conjugation. Hence a closure of an alternating word of
positive type is $+$amphichiral and one of negative type is $-$amphichiral.
\end{lemma}

\begin{proof}
The mirror image of $\sigma_1^{e_0}\sigma_2^{-e_1}\cdots$ is $\sigma_1^{-e_0}\sigma_2^{e_1}\cdots$; conjugation by $\Delta$ exchanges $\sigma_1$ and
$\sigma_2$, giving $\sigma_2^{-e_0}\sigma_1^{e_1}\sigma_2^{-e_2}\cdots$, a cyclic permutation of the alternating word of $\rot_1(s)$. Reading
backwards gives $\sigma_2^{-e_{L-1}}\sigma_1^{e_{L-2}}\cdots\sigma_1^{e_0}$, a cyclic permutation of the alternating word of $\rot_1(\rev s)$; its
closure is the knot with reversed orientation (rotate the closed braid by $\pi$ about a line that meets the braid axis at a right angle: the letters keep their signs and indices and their cyclic order is reversed). If
$s=\rot_r(s)$ with $r$ odd, then $\rot_1(s)=\rot_{1-r}(s)$ with $1-r$ even, so the mirror image is isotopic to the knot. If $\rev s$ is $s$ rotated by
an even amount, the mirror image and the reverse are closures of conjugate words.
\end{proof}

\begin{lemma}[negative type]\label{lem:N}
If $s$ is of negative type and the closure of its alternating word is a knot, it has $\det=(a-d)^2+(b+c)^2$ and $\det+4=(a+d)^2+(b-c)^2$ for a nonnegative integer
matrix $W=\begin{pmatrix}a&b\\c&d\end{pmatrix}$ with $ad-bc=1$.
\end{lemma}

\begin{proof}
Let $e_i=e_{k-i}$ with $k$ odd, $k\in\{1,\dots,L-1\}$, and $r=(k+1)/2$. Then $f=\rot_r(s)$ satisfies $f_i=f_{L-1-i}$, i.e.\ $f=(h,\rev h)$ with
$h=(f_0,\dots,f_{L/2-1})$. The cyclically rotated product of $M(s)$ starting at index $r$ equals $M(f)$ if $r$ is even and $J\,M(f)\,J$ if $r$ is odd
(since $JU(e)J=L(e)$); in both cases $\tr M(f)=\tr M(s)$. Let $W=\prod_{i<L/2}F_i(f_i)$ with $F_i=U$ for $i$ even and $F_i=L$ for $i$ odd. Since
$U(e)\T=L(e)$ and the indices $L/2-1-j$ and $L/2+j$ have different parity, $W\T$ is the product of the remaining factors of $M(f)$, so
$M(f)=WW\T$. $W$ is a product of nonnegative unimodular matrices, so $\det W=1$ and $\tr(WW\T)=a^2+b^2+c^2+d^2$. Hence
$\det=\tr M-2=a^2+b^2+c^2+d^2-2(ad-bc)=(a-d)^2+(b+c)^2$ and $\det+4=(a+d)^2+(b-c)^2$.
\end{proof}

\begin{lemma}[positive type]\label{lem:P}
If $s$ is of positive type and the closure of its alternating word is a knot, its determinant $n$ is a sum of two squares, and $n+4$ is a sum of two squares.
\end{lemma}

\begin{proof}
Let $s=\rot_r(s)$ with $r$ odd. The self-rotations of $s$ form a subgroup of $\Z/L$, so $g=\gcd(r,L)$ is a self-rotation; $g$ is odd and $L=gm$
with $m$ even. Let $V$ be the product of the first $g$ factors of $M(s)$. The $j$-th block of $g$ factors starts at index $jg$, which has the parity of
$j$, so the blocks are alternately $V$ and $JVJ$, and $M(s)=(V\cdot JVJ)^{m/2}=(VJ)^m$ with $\det(VJ)=-1$. Write $m=2^\nu o$ with $o$ odd, $\nu\ge1$,
$N_0=(VJ)^o$ ($\det N_0=-1$), $x=\tr N_0$, and $\tau_j=\tr N_0^{2^j}$. From $\tr(X^2)=(\tr X)^2-2\det X$: $\tau_1=x^2+2$ and $\tau_{j+1}=\tau_j^2-2$ for $j\ge1$.
Hence $n=\tau_\nu-2$ equals $x^2$ for $\nu=1$, equals $(x^2+2)^2-4=x^2(x^2+4)$ for $\nu=2$, and for $\nu\ge3$, by induction using
$\tau_{\nu}-2=(\tau_{\nu-1}-2)(\tau_{\nu-1}+2)$ and $\tau_{\nu-1}+2=\tau_{\nu-2}^2$,
\[
n=S^2(x^2+4),\qquad S=x\prod_{i=1}^{\nu-2}\tau_i \quad(\nu\ge2),
\]
so $n=(Sx)^2+(2S)^2$. Finally $n+4=x^2+2^2$ for $\nu=1$ and $n+4=\tau_\nu+2=\tau_{\nu-1}^2$ for $\nu\ge2$.
\end{proof}

\begin{lemma}[converse]\label{lem:C}
If $n$ is odd and $n$ and $n+4$ are sums of two squares, there is a sequence of negative type whose alternating word closes to a knot with
determinant~$n$.
\end{lemma}

\begin{proof}
Write $n=u^2+v^2$, $n+4=p^2+q^2$ with $u,v,p,q\ge0$. Each representation has one odd and one even component; let $u_o,p_o$ be the odd ones. If $p_o>u_o$, let $u,p$ be
the odd components (then $p\ge u+2$ and $q^2<v^2$); otherwise let $u,p$ be the even ones (then $q=p_o\le v=u_o$ and $p^2-u^2\ge4$). Put
$a=(p+u)/2$, $d=(p-u)/2$, $b=(v+q)/2$, $c=(v-q)/2$. These are nonnegative integers with $ad-bc=(p^2+q^2-u^2-v^2)/4=1$, $b\ge1$ (in both cases
$v\ge1$: in the first $v^2>q^2$, in the second $v$ is odd), so $W=\begin{pmatrix}a&b\\c&d\end{pmatrix}\ne I$, and $a-c=(p+u+q-v)/2>0$ (as $p+q\ge\sqrt{n+4}>v$).

\emph{Row dominance.} A nonnegative integer matrix $W\ne I$ with $\det W=1$ has one row entrywise $\ge$ the other: otherwise either $a>c$, $b<d$ and
$ad-bc\ge(c+1)(b+1)-bc=b+c+1$, which equals $1$ only for $W=I$; or $a<c$, $b>d$ and $ad-bc\le(c-1)(b-1)-bc\le0$; and both rows dominating each other
gives $\det W=0$. If the top row dominates, $W=U(1)W'$, otherwise $W=L(1)W'$, with $W'$ nonnegative, $\det W'=1$ and entry sum smaller by $c+d$ resp.\
$a+b$, which is $\ge1$; so the process terminates at $W'=I$, and $W$ is a product of factors $U(1)$ and $L(1)$. Collecting equal neighbours,
$W=U(h_0)L(h_1)\cdots$ with all $h_i\ge1$; it starts with $U$ because $a>c$. Then $s=(h,\rev h)$ is of negative type with $M(s)=WW\T$ as in
Lemma~\ref{lem:N}, so its determinant is $(a-d)^2+(b+c)^2=u^2+v^2=n$.
\emph{Knot.} $\tr M(s)=n+2$ is odd, so $M(s)\bmod2$ has order $3$ in $SL(2,\mathbb F_2)\cong S_3$; since $\sigma_1\mapsto U(1)$ and
$\sigma_2^{-1}\mapsto L(1)$ modulo $2$ send the generators to two distinct transpositions, this is the permutation of the braid up to an automorphism
of $S_3$, and it is a $3$-cycle.
\end{proof}

\section{Proof of Theorem \ref{thm:main}}

\begin{proposition}[{\cite[Thm.~1.1(a)]{St15}}; for prime knots {\cite[Theorems~4.13 and B]{J}}]\label{prop:A}
An amphichiral knot $K$ of braid index $3$ is the closure of an alternating word of positive or negative type.
\end{proposition}

\begin{proof}[Derivation from~\cite{BM}]
Let $W$ be the Schreier normal form of a 3-braid representing $K$. The closures of $C^k\sigma_1^p$ and $C^k\sigma_1\sigma_2\sigma_1$ are not knots
(the permutation is not a $3$-cycle), and $C^k\sigma_1\sigma_2$, $C^k\sigma_1\sigma_2\sigma_1\sigma_2$ close to torus knots $T(3,m)$, $3\nmid m$, which
have braid index $\le2$ for $|m|\le2$; for $|m|\ge4$ they cannot be amphichiral of braid index $3$: $W$ would be a minimal braid with exponent sum
$6k+2$ resp.\ $6k+4\ne0$, while a minimal braid of an amphichiral knot has exponent sum $0$ (by Dynnikov--Prasolov~\cite[Theorem~9]{DP}, as in~\cite[Lemma~3.3 and eq.~(1)]{J}). (They are
chiral in any case; their signatures were computed by Hirzebruch, Murasugi and Shinohara, as cited in~\cite[\S2]{E}.) So $W=C^k\sigma_1^{-p_1}\sigma_2^{q_1}\cdots\sigma_1^{-p_r}\sigma_2^{q_r}$ with array
$A=(p_1,q_1,\dots,p_r,q_r)$. Let $W^*$ be the mirror word and $W^{\leftarrow}$ the word read backwards. If $K$ is $+$amphichiral, $W^*$ and $W$ represent
the same oriented knot; if $-$amphichiral, $W^*$ and $W^{\leftarrow}$ do. By the classification~\cite[pp.~76--77, 106--107]{BM} (braid index $3$; in the
flype case a word and its reverse represent the same knot) this gives $W^*\sim W$ or $W^*\sim W^{\leftarrow}$ in both cases. The normal form of $W^*$ is
$C^{-k}$ with array $\rot_1(A)$~\cite[p.~107]{BM}, that of $W^{\leftarrow}$ is $C^k$ with array $\rot_1(\rev A)$ (reading backwards gives $C^k\sigma_2^{q_r}\sigma_1^{-p_r}\cdots\sigma_2^{q_1}\sigma_1^{-p_1}$, cyclically $C^k\sigma_1^{-p_r}\sigma_2^{q_{r-1}}\cdots\sigma_1^{-p_1}\sigma_2^{q_r}$), and two generic normal forms are
conjugate iff the powers of $C$ agree and the arrays agree up to an even rotation. Hence $k=0$, and $A=\rot_{\mathrm{odd}}(A)$ or $A$ is $\rev A$ rotated
by an even amount. With $k=0$, $W$ is the mirror of the alternating word of $A$, i.e.\ (Lemma~\ref{lem:S}) conjugate to the alternating word of
$\rot_1(A)$, and both conditions are invariant under rotation.
\end{proof}

\begin{proof}[Proof of Theorem~\ref{thm:main}]
``$\subseteq$'': braid index $1$ gives $n=1=1^2+0^2$, $n+4=5=1^2+2^2$; braid index $2$ gives no amphichiral knot; braid index $3$ is covered by
Proposition~\ref{prop:A} with Lemmas~\ref{lem:N} and~\ref{lem:P}. ``$\supseteq$'': Lemma~\ref{lem:C} gives a word of negative type closing to a knot
with determinant $n$; the knot is $-$amphichiral by Lemma~\ref{lem:S} and has braid index $\le3$.
\end{proof}

\section{Remarks}

\begin{remark}[composite knots and small determinants]
The closure of $\sigma_1^p\sigma_2^{-p}$ ($s=(p,p)$, $p$ odd) is the connected sum of $T(2,p)$ with its mirror image, with determinant $p^2$. For
$n=1,9,49$ --- the three values excluded in Stoimenow's conjecture on prime alternating achiral knots~\cite{StCAG} --- the only sequences of positive or
negative type, up to rotation, are $(1,1)$, $(3,3)$, $(7,7)$: the unknot, and the connected sums of $3_1$ resp.\ $7_1$ with their mirror images. The search
is finite because $\det\ge(\sum_{i\text{ even}}e_i)(\sum_{i\text{ odd}}e_i)$: the trace is a polynomial in the exponents with nonnegative coefficients and constant term $2$ (all exponents $0$ give the identity), in which every product of one $U$- and one $L$-exponent occurs with
coefficient $1$.
\end{remark}

\begin{remark}[Stoimenow's conjecture]\label{rem:stoimenow}
Stoimenow conjectured~\cite[Conjecture~4.20; Conjecture~4.2 in arXiv v3]{StCAG}: an odd $n$ is the determinant of a prime alternating achiral knot iff $n$
is a sum of two squares and $n\notin\{1,9,49\}$; he settled the $n$ with a representation by coprime components and the squares, and reports the
conjecture true for $n\le2000$. Up to $10{,}693$ (the largest determinant among the verified amphichiral prime alternating knots of Burton's census,
which go up to $18$ crossings --- alternating amphichiral knots have even crossing number) there remain $169$ odd $n>2000$ not settled by him (sums of two squares that are not squares and have no representation by coprime components); for $74$ of them $n+4$ is
not a sum of two squares, so by Corollary~\ref{cor:b5} any witness has braid index at least $5$. Ferudun~\cite{F} realizes every $n\le10^6$ left open
by Stoimenow (unrefereed). Of the $74$, $26$ have no witness in the census: $6309$, $6413$, $7281$, $7301$, $7857$, $8361$, $8433$, $8477$, $8793$, $8973$, $8993$, $9065$, $9081$, $9117$,
$9225$, $9245$, $9297$, $9477$, $9621$, $9873$, $10125$, $10161$, $10377$, $10413$, $10469$, $10629$. Forty-one further alternating census knots with $14$ to $18$
crossings, amphichiral according to SnapPy's numerical symmetry group, stayed undecided in our verified computation; none of their determinants is
among these $26$. For $25$ of the $26$ we found witnesses that are closures of $5$-braid words $U\rho(U)$, where $U$ is a word in
$\sigma_1,\sigma_2^{-1},\sigma_3,\sigma_4^{-1}$ and $\rho(U)$ is $U$ read backwards with each $\sigma_i^{\pm1}$ replaced by $\sigma_{5-i}^{\mp1}$; they are prime
(verified hyperbolic), alternating (reduced alternating diagrams of $20$ to $26$ crossings), verified $-$amphichiral, with Morton--Franks--Williams
bound $5$, so of braid index exactly $5$; details in the deposit~\cite{Z}. The remaining value, $n=6413$, is covered only by~\cite{F}. No knot that is the closure of a word $U\rho(U)$ with $U$ of at most $20$ letters (all $4^m$ words $U$ for each $1\le m\le20$; their closures
are braid diagrams of $2m\le40$ crossings) has determinant $6413$.
\end{remark}

\begin{remark}[the alternating half without the classification]\label{rem:erle}
That an amphichiral knot of braid index $3$ is alternating also follows without the classification theorem of~\cite{BM}. The non-generic normal forms
are excluded as in the first lines of the proof of Proposition~\ref{prop:A}. For a generic form $W=C^k\sigma_1^{-p_1}\sigma_2^{q_1}\cdots\sigma_1^{-p_r}
\sigma_2^{q_r}$ write $P=\sum p_i$, $Q=\sum q_i$; Erle's Theorem~2.6(iv)~\cite{E} states $\sigma=P-Q-4k$, which equals $2k-e(W)$, $e$ the exponent sum; his printed proof treats $Q\ge3$ with $k\ge1$
and with $k\le-1$ (the remaining cases are left as an exercise); we use only the case $k\le-1$, $Q\ge3$, and only that $\sigma\ne0$ there, which does not depend on the sign convention for $\sigma$. For a
minimal braid of an amphichiral knot $e(W)=0$ (by Dynnikov--Prasolov~\cite[Theorem~9]{DP}, as in~\cite[Theorem~3.2, Lemma~3.3, eq.~(1)]{J}) and $\sigma=0$. Suppose $k\ne0$; then
$P=Q+6k$. If $k\le-1$,
$Q=P-6k\ge7$; if $k\ge1$, the mirror of $W$, conjugated by $\Delta$ and cyclically rotated, is the generic form with $C^{-k}$ and the roles of the $p_i$ and
$q_i$ exchanged, whose $Q$ is $P\ge7$, and $\sigma$ changes sign under mirroring. Either way $\sigma=\pm2k\ne0$, a contradiction; so $k=0$ and $W$ is
alternating. The same argument with the extreme degrees of the Jones polynomial in place of Dynnikov--Prasolov is Chbili's~\cite[Proposition~5.1]{Ch}.
Stoimenov's proof of~\cite[Theorem~1.1]{St15} also bypasses the classification, as described in \S2.
\end{remark}

\begin{remark}[context]
Datta~\cite{D} relates links of braid index $\le3$ to binary quadratic forms of discriminant $\theta^2-4$, $\theta$ the trace of the Burau matrix at $t=-1$;
for the knots here $\theta=n+2$ and the discriminant is $n(n+4)$. He does not consider amphichirality. For $\nu=1$ in Lemma~\ref{lem:P} the
alternating word is of the form $\gamma\bar\gamma$, with $\bar\gamma$ obtained from $\gamma$ by $\sigma_i^{\pm1}\mapsto\sigma_{3-i}^{\mp1}$ (as in~\cite[p.~10]{StW}),
and the square determinant can be explained as Stoimenow does for the Turk's head braids $(\sigma_1\sigma_2^{-1})^k$ with $k$ odd (this is the case $\nu=1$; the closure is a knot iff $3\nmid k$), through strong achirality and Hartley--Kawauchi~\cite[Theorems~3.1, 3.2]{StW}.
\end{remark}

\section{Checks}\label{sec:checks}

\emph{Braid index test of ``$\subseteq$''.} Contrapositive: an amphichiral knot whose determinant $n$ has $n+4$ not a sum of two squares must have
braid index $\ge5$. (a) KnotInfo~\cite{KI}, $\le12$ crossings, tabulated braid index: of the 78 amphichiral knots, the 22 of braid index $\le3$ all
satisfy the condition, and the 22 that violate it (a different set of knots) all have braid index $\ge5$; without amphichirality the implication fails ($107$ chiral knots of braid
index $3$ have $n+4$ not a sum of two squares, the first being $5_2$). (b) Burton's census~\cite{Bu}, $\le19$ crossings, amphichirality verified with
SnapPy~\cite{SP} (knots left undecided by the verified computation are not counted): $3{,}049$ alternating and $1{,}950$ non-alternating knots violate the condition. The Morton--Franks--Williams lower bound is $\ge4$ for all
alternating ones and for $1{,}892$ of the others; the remaining $58$ have Jones span smaller than the crossing number, so they are not alternating (the
census crossing numbers are exact), and by the argument of Remark~\ref{rem:erle}, which does not use the classification of~\cite{BM}, they do not have braid index~$3$. None of the
$4{,}999$ has braid index~$3$; with the odd braid index~\cite[Theorem~A]{J}, all have braid index $\ge5$. These checks are consistent with Theorem~\ref{thm:main}; by their nature they cannot prove it.

\emph{Other checks.} Lemmas~\ref{lem:N}, \ref{lem:P}, \ref{lem:C} and the type conditions were checked by exhaustive enumeration on small boxes, and the
relation $\det=|2-\tr M|$ against SnapPy for random 3-braid knots; scripts and logs are in the deposit~\cite{Z}.

\subsection*{Tool use}
This work was prepared with substantial use of Claude Code (Anthropic), an agentic AI system based on large language models, working together with
the author.
\begin{itemize}[leftmargin=1.4em]
\item \emph{Mathematics.} The proofs were worked out by Claude Code under the author's direction.
\item \emph{Writing and typesetting.} The text and the \LaTeX{} source were written by Claude Code, in English.
\item \emph{Code.} The checks were written and run by Claude Code in Python, with SnapPy~\cite{SP}, Sage and the KnotInfo database~\cite{KI}; the
$5$-braid search of Remark~\ref{rem:stoimenow} ran as a CUDA kernel through CuPy on a GPU.
\item \emph{Checking.} The argument and drafts of this note were read by separate Claude Code agents in independent cold reads and by outside AI readers (Claude in a web
browser, DeepSeek and ChatGPT). Their findings were checked and incorporated where they held.
\item \emph{Literature.} The search was done by a Claude Code agent (web and bibliographic databases), with the author by hand (Google Scholar,
MathOverflow) and, by Claude Code, through the arXiv API; an outside AI reader (Claude in a web browser) suggested further sources. \cite{StW} was found in the author's search and \cite{St15}
was pointed out by an outside AI reader (Claude in a web browser), both only after the argument was written.
\end{itemize}
The author takes full responsibility for all contents.

\subsection*{Data and code}
Scripts and their outputs are in the deposit~\cite{Z} on Zenodo; its README describes which files support the computational checks.


\begin{thebibliography}{StCAG}

\bibitem[BM]{BM} J.~S. Birman, W.~W. Menasco, \emph{Studying links via closed braids III: Classifying links which are closed 3-braids}, Pacific J.\ Math.\
\textbf{161} (1993) 25--113.

\bibitem[Bu]{Bu} B.~A. Burton, census of prime knots up to 19 crossings, data files \url{https://people.smp.uq.edu.au/BenjaminBurton/knots/}, linked
from \url{https://regina-normal.github.io/data.html} (downloaded 2026).

\bibitem[Ch]{Ch} N.~Chbili, \emph{A note on the Jones polynomials of 3-braid links}, Siberian Math.\ J.\ \textbf{63} (2022) 983--994.

\bibitem[D]{D} A.~Datta, \emph{A novel connection between integral binary quadratic forms and knot polynomials}, arXiv:2204.13660v1 (2022).

\bibitem[DP]{DP} I.~Dynnikov, M.~Prasolov, \emph{Bypasses for rectangular diagrams. Proof of Jones' conjecture and related questions}, arXiv:1206.0898v2 (2013);
journal reference given there: Proc.\ of Moscow Math.\ Soc.\ \textbf{74} (2013), no.~1, 115--173; English translation: Trans.\ Moscow Math.\ Soc.\ \textbf{74} (2013) 97--144.

\bibitem[E]{E} D.~Erle, \emph{Calculation of the signature of a 3-braid link}, Kobe J.\ Math.\ \textbf{16} (1999) 161--175.

\bibitem[F]{F} A.~Ferudun, \emph{Determinants of prime alternating achiral knots: a partial answer to a problem of Stoimenow}, Zenodo,
\href{https://doi.org/10.5281/zenodo.23482759}{10.5281/zenodo.23482759} (2026).

\bibitem[HK]{HK} R.~Hartley, A.~Kawauchi, \emph{Polynomials of amphicheiral knots}, Math.\ Ann.\ \textbf{243} (1979) 63--70.

\bibitem[J]{J} Hyungseok Jung, \emph{Amphichiral knots: odd braid index and symmetry classification in the three-braid case}, arXiv:2609.12295v1 (2026).

\bibitem[KI]{KI} C.~Livingston, A.~H. Moore, \emph{KnotInfo: Table of Knot Invariants}, knotinfo.org (accessed 2026).

\bibitem[Li]{Li} P.~Lisca, \emph{On 3-braid knots of finite concordance order}, Trans.\ Amer.\ Math.\ Soc.\ \textbf{369} (2017), no.~7, 5087--5112; arXiv:1504.03129v3.

\bibitem[SP]{SP} M.~Culler, N.~M. Dunfield, M.~Goerner, J.~R. Weeks, \emph{SnapPy}, version 3.3.2.

\bibitem[St15]{St15} A.~Stoimenov, \emph{The classification of partially symmetric 3-braid links}, Open Math.\ \textbf{13} (2015) 444--470,
\href{https://doi.org/10.1515/math-2015-0043}{10.1515/math-2015-0043}.

\bibitem[StCAG]{StCAG} A.~Stoimenow, \emph{Square numbers, spanning trees and invariants of achiral knots}, Comm.\ Anal.\ Geom.\ \textbf{13} (2005) 591--631;
Conjecture~4.20 there; also arXiv:math/0003172v3, Conjecture~4.2 there.

\bibitem[StMZ]{StMZ} A.~Stoimenow, \emph{Square numbers and polynomial invariants of achiral knots}, Math.\ Z.\ \textbf{255} (2007) 703--719.

\bibitem[St06]{St06} A.~Stoimenow, \emph{Non-triviality of the Jones polynomial and the crossing numbers of amphicheiral knots}, arXiv:math/0606255v2 (2007).

\bibitem[StW]{StW} A.~Stoimenow, \emph{Wheel graphs, Lucas numbers and the determinant of a knot}, preprint, Max-Planck-Institut f\"ur Mathematik (2000).

\bibitem[Z]{Z} W.~A. Stanggassinger, \emph{Determinants of amphichiral knots of braid index at most three}, note, scripts, logs and result files, Zenodo,
\href{https://doi.org/10.5281/zenodo.24157621}{10.5281/zenodo.24157621} (2026).

\end{thebibliography}
\end{document}